\documentclass[12pt]{amsart}
\usepackage{amscd,amssymb,amsthm,verbatim,slashed}
\newcommand{\Area}{\operatorname{Area}}

\newcommand{\Diff}{\operatorname{Diff}}

\newcommand{\dist}{\operatorname{dist}}

\newcommand{\dvol}{\operatorname{dvol}}

\newcommand{\HH}{\operatorname{H}}

\newcommand{\Nil}{\operatorname{Nil}}

\newcommand{\OO}{\operatorname{O}}

\newcommand{\R}{{\mathbb R}}

\newcommand{\SL}{\operatorname{SL}}
\newcommand{\SO}{\operatorname{SO}}
\newcommand{\Sol}{\operatorname{Sol}}

\newcommand{\Z}{{\mathbb Z}}

\numberwithin{equation}{section}
\theoremstyle{plain}

\newtheorem{theorem}{Theorem}

\theoremstyle{remark}

\newtheorem{remark}{Remark}

\usepackage{graphicx}
\input{amssym.def}
\input{amssym}
\input{epsf}
\usepackage{a4wide}

\begin{document}

\title[Remark about scalar curvature on certain noncompact manifolds]
      {Remark about scalar curvature on certain noncompact manifolds}

\author{John Lott}
\address{Department of Mathematics\\
University of California, Berkeley\\
Berkeley, CA  94720-3840\\
USA} \email{lott@berkeley.edu}

\begin{abstract}
We give a sufficient condition to rule out complete Riemannian metrics with
nonnegative scalar curvature on the interiors of handlebodies.  In higher dimensions,
we give examples of ends of manifolds with positive scalar curvature metrics.
\end{abstract}

\date{April 27, 2026}

\maketitle

\section{Introduction} \label{sect1}

Let $\overline{M}$ be an orientable three dimensional handlebody, i.e. a repeated boundary connected sum of solid tori. 
Equivalently, $\overline{M}$ is homeomorphic to a regular neighborhood of a bouquet of circles in $\R^3$. 
We say that $\overline{M}$ has genus $g$ if
$\partial \overline{M}$ is a surface of genus $g$. 

Let $M$ be the interior of $\overline{M}$. If $g = 0$ then $M$
is diffeomorphic to $\R^3$ and has a complete Riemannian metric of uniformly positive scalar curvature.  If $g=1$ then $M$ is diffeomorphic to $S^1 \times \R^2$ and has a complete Riemannian metric of nonuniformly positive scalar curvature.
If $g > 1$ then it is unclear whether $M$ admits a complete Riemannian metric of nonnegative scalar curvature.

Chodosh, Lai and Xu used the inverse mean curvature flow to
show that this is not the case under the condition that the metric has bounded sectional curvature and positive injectivity radius \cite{Chodosh-Lai-Xu (2025)}.  Their proof actually shows that the metric cannot have 
nonnegative scalar curvature outside of a compact set,
indicating that there is an end obstruction.  We give another condition to rule out nonnegative scalar curvature
outside of a compact subset.  The condition roughly says that there are surfaces going out the end whose areas
do not increase too fast.

\begin{theorem} \label{thm1}
Let $\overline{M}$ be an orientable three dimensional handlebody. Let $M$ be the interior of $\overline{M}$.  Suppose that $M$ has a complete Riemannian metric with nonnegative scalar curvature outside of a compact subset.  Choose a basepoint $m_0 \in M$. Given $r > 0$, let $A(r)$ be the infimum of the areas of smooth closed connected embedded surfaces in $M$
which are homologous to $\partial \overline{M}$ in a collar end $[0,1)\times \partial\overline{M}$
and lie outside of $B(m_0, r)$.
Suppose that
$\liminf_{r \rightarrow \infty} r^{-2} A(r) < \frac{12}{\pi}$.  Then the genus of $\overline{M}$ is at most one.
\end{theorem}

The proof of Theorem \ref{thm1} indicates that a possible end obstruction depends on whether the boundary at infinity has almost nonnegative scalar curvature, in a suitable sense.  
The next result, which is valid in any dimension, gives a
sufficient condition to have positive scalar curvature on an end.

\begin{theorem} \label{thm2}
Let $N$ be a compact manifold that admits a Riemannian metric with nonnegative scalar curvature.  Let $X$ be the total space of an $S^1$-bundle over $N$. Then $[0, \infty) \times X$ admits a complete Riemannian metric with positive scalar curvature.
\end{theorem}

The organization of the paper is as follows.  In Section \ref{sect2} we prove Theorem \ref{thm1}.  
We also state a band width inequality.
In Section \ref{sect3} we
prove Theorem \ref{thm2}. Section \ref{sect4} has a discussion.

I thank Otis Chodosh, Yi Lai and Kai Xu for conversations.

\section{Proof of Theorem \ref{thm1}} \label{sect2}

The proof uses $\mu$-bubbles, for which references are \cite[Section 3]{CL (2024)}
and \cite[Section 5]{Gromov (2023)}.
We find quantitative negative lower bounds on the total
curvatures for a sequence of $\mu$-bubbles exiting the end.  
One ingredient is upper area bounds for the $\mu$-bubbles.

By assumption, there are some $\epsilon > 0$ and a sequence $r_i \rightarrow \infty$ so that for each $i$, there is
a smooth closed connected embedded surface $\Sigma_i \subset M$ with
\begin{itemize}
\item $\Sigma_i \subset M - B(m_0, r_i)$,
\item $\Sigma_i$ is homologous to $\partial \overline{M}$ in the collar end 
$[0,1)\times \partial\overline{M}$, and
\item $\Area(\Sigma_i) \le \left( \frac{12}{\pi} - \epsilon \right) r_i^2$.
\end{itemize}
Let $\epsilon^\prime > 0$ be a small parameter that will be independent of $i$.
For the moment, we drop the index $i$. 
Choose $\epsilon^{\prime \prime} > 0$, a small parameter (which can depend on $i$).
Put $D = \max_{p \in \Sigma} d(m_0, p)$.
The surface $\Sigma$ is separating.
Let $d_\Sigma$ denote the signed distance from $\Sigma$, with $d_\Sigma(m)$
going to $- \infty$ as $m$ exits the end of $M$. There is some $\delta > 0$ so that $d_\Sigma$ is smooth in the
$\delta$-neighborhood of $\Sigma$. Furthermore, we can find some $\widehat{d}_\Sigma \in C^\infty(M)$ so that
on $B(m_0, 3D)$,
\begin{itemize}
\item $\widehat{d}_\Sigma^{-1}(0) = \Sigma$,
\item $|\widehat{d}_\Sigma - d_\Sigma| \le \epsilon^{\prime \prime}$,
\item $|\nabla \widehat{d}_\Sigma| \le 1 + \epsilon^{\prime \prime}$, and
\item $\widehat{d}_\Sigma$ equals $d_\Sigma$ on the $\frac{\delta}{2}$-neighborhood of $\Sigma$.
\end{itemize}

By Sard's theorem, we can find $L \in [(1 - \epsilon^\prime) r - \epsilon^{\prime \prime}, (1 - \epsilon^\prime) r]$ so that
$-L$ and $L$ are regular values for $\widehat{d}_\Sigma$. Put
\begin{itemize}
\item $N$ to be the connected component of  $\widehat{d}_\Sigma^{-1} [-L, L]$ containing $\Sigma$,
\item $\partial_+ N = \widehat{d}_\Sigma^{-1}(L) \cap N$ and $\partial_- N = \widehat{d}_\Sigma^{-1}(-L) \cap N$, and
\item $\Omega_0 = \widehat{d}_\Sigma^{-1} [0, L] \cap N$.
\end{itemize}
The boundary component $\partial_- N$ is farther out the end than $\partial_+ N$.

If $i$ is large and $\epsilon^{\prime \prime}$ is small then $N$ has nonnegative scalar curvature.
Put $h = \frac{(1 + \epsilon^{\prime \prime}) 2\pi}{3L} \tan \left( \frac{\pi}{2L} \widehat{d}_\Sigma \right)$ on the interior of $N$.
Then $h(m)$ approaches $\pm \infty$ as $m$ approaches $\partial_\pm N$. 

Let $\Omega$ be a codimension-zero submanifold-with-boundary of $N$ so that the symmetric difference
$\Omega \Delta \Omega_0$ has compact support in the interior of $N$. In particular, $\partial_+ N \subset
\partial \Omega$. Put $\partial_* \Omega = \partial \Omega - \partial_+ N$. Define the functional
\begin{equation} \label{2.1}
{\mathcal A}(\Omega) = A(\partial_* \Omega) - \int_N (\chi_\Omega - \chi_{\Omega_0}) h \dvol_N.
\end{equation}
By the standard $\mu$-bubble existence theorem with blow-up barriers at $\partial_\pm N$,
there is a smooth minimizer $\Omega$ for ${\mathcal A}$. The mean curvature of $\partial_* \Omega$ is
$h \Big|_{\partial_* \Omega}$. The second variation formula for ${\mathcal A}$ implies that for any $\psi \in C^\infty(\partial_* \Omega)$, we have
\begin{equation} \label{2.2}
0 \le \int_{\partial_* \Omega} \left( |\nabla_{\partial_* \Omega} \psi|^2 - \frac12
\left( R_M - R_{\partial_* \Omega} + \frac32 h^2 + 2 \langle \nabla h, \nu \rangle \right) \psi^2 \right) dA,
\end{equation}
where $\nu$ is the outward pointing normal
vector.  In particular, if $C$ is a connected component of $\partial_* \Omega$ then taking $\psi = \chi_C$ gives 
\begin{equation} \label{2.3}
\int_C R_C \: dA \ge \int_C \left( \frac32 h^2 - 2 |\nabla h| \right) \: dA.
\end{equation}

Now
\begin{align} \label{2.4}
\frac32 h^2 - 2 |\nabla h| & = \frac{2(1 + \epsilon^{\prime \prime})^2 \pi^2}{3L^2} \tan^2 \left( \frac{\pi}{2L} \widehat{d}_\Sigma \right) 
- \frac{2(1 + \epsilon^{\prime \prime}) \pi^2}{3L^2}  |\nabla 
\widehat{d}_\Sigma | \sec^2 \left( \frac{\pi}{2L} \widehat{d}_\Sigma \right) \\
& \ge \frac{2(1 + \epsilon^{\prime \prime})^2 \pi^2}{3L^2} \tan^2 \left( \frac{\pi}{2L} \widehat{d}_\Sigma \right) 
- \frac{2(1 + \epsilon^{\prime \prime}) \pi^2}{3L^2} (1 + \epsilon^{\prime \prime}) \sec^2 \left( \frac{\pi}{2L} \widehat{d}_\Sigma \right) \notag \\
& = - \: \frac{2(1 + \epsilon^{\prime \prime})^2 \pi^2}{3L^2}. \notag
\end{align}
Finally,
${\mathcal A}(\Omega) \le {\mathcal A}(\Omega_0) = A(\Sigma)$, so
\begin{align} \label{2.5}
A(\partial_* \Omega) & \le A(\Sigma) + \int_N (\chi_\Omega - \chi_{\Omega_0}) h \dvol_N \\
& = A(\Sigma) + \int_{\Omega - (\Omega_0 \cap \Omega)} h \dvol_N -  \int_{\Omega_0 - 
(\Omega_0 \cap \Omega)} h \dvol_N. \notag
\end{align}
If $\epsilon^{\prime \prime} \ll \delta$ then $h$ is nonpositive on $\Omega - (\Omega_0 \cap \Omega)$ and
nonnegative on $\Omega_0 - (\Omega_0 \cap \Omega)$, so $A(\partial_* \Omega) \le A(\Sigma)$. In particular,
$A(C) \le A(\Sigma)$.

Equations (\ref{2.3}) and (\ref{2.4}) now imply that $\int_C R_C \:  dA \ge - \: \frac{2(1 + \epsilon^{\prime \prime})^2 \pi^2}{3L^2} A(\Sigma)$.
Taking $\epsilon^{\prime \prime}$ sufficiently small ensures that $\int_C R_C \:  dA \ge - \: \frac{2\pi^2}{3(1 - \epsilon^\prime)^2 r^2} A(\Sigma)$.

Restoring the index $i$ shows that for $i$ large, any connected component $C$ of $\partial_* \Omega_i$ satisfies
\begin{equation} \label{2.6}
\int_{C} R_{C} \:  dA \ge - \: \frac{2\pi^2}{3(1 - \epsilon^\prime)^2 r_i^2} A(\Sigma_i)
\ge - \: \frac{2\pi^2}{3(1 - \epsilon^\prime)^2} \left( \frac{12}{\pi} - \epsilon \right).
\end{equation}
If $\epsilon^\prime$ is sufficiently small then $\int_{C} R_{C} \:  dA 
> - \: \frac23 \pi^2 \cdot \frac{12}{\pi} = - 8 \pi$. By the Gauss-Bonnet theorem,
$\int_{C} R_{C} \:  dA = 4 \pi \chi(C)$. Hence $C$ is diffeomorphic to $S^2$ or $T^2$.

For fixed $\epsilon^\prime > 0$,
there is a subset $S$ of $\overline{M}$, diffeomorphic to $[0,1] \times \partial \overline{M}$
and containing $\partial \overline{M}$, so that
$\partial_* \Omega_i$ lies in $S$ for all large $i$. Let $e_i : \partial_* \Omega_i \rightarrow S$ be the embedding and let
$\alpha : S \rightarrow \partial \overline{M}$ be the projection map.  Since the oriented surfaces
$\Sigma_i$ and $\partial_*\Omega_i$ are homologous in the  collar end,
the image of the total fundamental class
$[\partial_*\Omega_i]\in \HH_2(\partial_*\Omega_i; \Z)$
under $(\alpha\circ e_i)_*$ is $[\partial \overline{M}]$.
Equivalently, if $\{C_{i,j}\}$ are the connected components of $\partial_* \Omega_i$ and
$d_{i,j}$ is the degree of $\alpha\circ e_i|_{C_{i,j}}$, then
$\sum_j d_{i,j}=1$.
Hence some component $C_{i,j}$ has a map to 
$\partial \overline{M}$ of nonzero degree.  As $C_{i,j}$ is diffeomorphic to $S^2$ or $T^2$,
the genus of $\partial \overline{M}$ is at most one. This proves the theorem.

\begin{remark}
The $n$-dimensional analog is the following. Put 
$h = \frac{(1 + \epsilon^{\prime \prime}) (n-1)\pi}{nL} \tan \left( \frac{\pi}{2L} \widehat{d}_\Sigma \right)$.
Let $\psi$ be the eigenfunction corresponding to the lowest eigenvalue $\lambda$ of the operator
$- \triangle + \frac12 (R_\Sigma - R_M - \frac{n}{n-1} h^2 + 2 |\nabla h|)$. Then
$\lambda$ is nonnegative and $\psi$ is positive. Put
$R_1 = R_\Sigma - 2 \frac{\triangle \psi}{\psi}$; it is the modified scalar curvature $R_q$ associated to
the pair $(g_\Sigma, \psi \dvol_\Sigma)$ when $q = 1$, in the sense of \cite{Lott (2007)}. Then
$R_1 \ge - \: \frac{(n-1) (1+ \epsilon^{\prime \prime})^2 \pi^2}{nL^2}$.
\end{remark}

\begin{remark}
The proof of Theorem \ref{thm1} gives the following band width inequality.  Suppose that $M$ is diffeomorphic to
$[-1,1] \times V$, where $V$ is a compact orientable surface of genus $g \ge 1$. Put $\Sigma_x = \{x\} \times V$. 
Suppose that $R_M \ge r_0 > - \:
\frac{8 \pi (g-1)}{A(\Sigma_0)}$. Then
\begin{equation}
\min(\dist(\Sigma_0, \Sigma_{-1}), \dist(\Sigma_0, \Sigma_1)) \le  \pi  \sqrt{
\frac{
2A(\Sigma_0)
}{
24 \pi (g-1) + 3r_0 A(\Sigma_0)
}
}
.
\end{equation}
One sees that if $g > 1$ then this has a weaker assumption and stronger conclusion than the $g=1$ case, 
the latter of which is
area independent.

If instead $M$ is diffeomorphic to $[0,1] \times V$, and $\Sigma_0$ has nonnegative mean curvature, then a
doubling argument gives the same inequality for
$\dist(\Sigma_0, \Sigma_1)$.
\end{remark}

\section{Proof of Theorem \ref{thm2}} \label{sect3}

Let $N^{\,n-2}$ be compact with a fixed Riemannian metric $h$ and scalar curvature $R_h\ge 0$.
We can assume that $N$ is connected.
Let $\pi:X\to N$ be an $S^1$-bundle. We can reduce the structure group of $X$ from $\Diff(S^1)$
to $\OO(2)$. After passing to a double cover if necessary and working equivariantly, we can
assume that the structure group is $\SO(2)$. Then $X$ is a principal $\SO(2)$-bundle over $N$.
Let $\theta$ be a connection $1$-form on $X$, with curvature
\begin{equation} \label{3.1}
d\theta=\pi^*\Omega,\qquad \Omega\in\Omega^2(N).
\end{equation}
Put $\Omega_{ij} = \Omega(\partial_i, \partial_j)$ and $|\Omega|_h^2 = \frac12 h^{ik} h^{jl} \Omega_{ij} \Omega_{kl}$.
Define
\begin{equation} \label{3.2}
\Omega_\infty=\max_{N}|\Omega|_h<\infty .
\end{equation}
For functions $a(t)>0$ and $b(t)>0$ on $[0,\infty)$, consider on $[0,\infty)\times X$ the metric
\begin{equation} \label{3.3}
\bar g \;=\; dt^2 + g(t),\qquad 
g(t)=a(t)^2\,\theta^2+b(t)^2\,\pi^*h .
\end{equation}
All norms $|\cdot|_h$ below are computed using $h$, and primes denote $t$-derivatives.
The formulas below were generated with artificial intelligence and checked with natural perseverance.

The scalar curvature of $\bar g$ on $[0,\infty)\times X$ is
\begin{equation} \label{3.4}
\begin{aligned}
R_{\bar g}
&=\frac{1}{b^2}R_h-\frac{a^2}{4b^4}|\Omega|_h^2
-2\frac{a''}{a}-2(n-2)\frac{b''}{b}
-2(n-2)\frac{a'}{a}\frac{b'}{b}
-(n-2)(n-3)\Big(\frac{b'}{b}\Big)^2.
\end{aligned}
\end{equation}
Since $R_h\ge 0$, for lower bounds one may discard the nonnegative term $\frac{1}{b^2}R_h$.

We now separate into the cases $n=2$, $n=3$, $n=4$, $n=5$ and $n \ge 6$.

\subsection*{${\bf n=2}$}
Here $\dim N=0$, so $R_h\equiv0$ and $\Omega\equiv0$. The metric is $\bar g=dt^2+a(t)^2\,d\phi^2$ and
\begin{equation} \label{3.5}
\;R_{\bar g}(t)=-2\,\frac{a''(t)}{a(t)}.\;
\end{equation}
Example:
\begin{equation} \label{3.6}
a(t)=1+\sqrt{1+t}\quad\Longrightarrow\quad
R_{\bar g}(t)=\frac{1}{2(1+t)^{3/2}\big(1+\sqrt{1+t}\big)}\;>\;0.
\;
\end{equation}

\subsection*{${\bf n=3}$}
Now $\dim N=1$, so $R_h\equiv0$ and $\Omega\equiv0$.
Taking $b(t)\equiv 1$ gives again
\begin{equation} \label{3.7}
\;R_{\bar g}(t)=-2\,\frac{a''(t)}{a(t)}.\;
\end{equation}
Equation (\ref{3.6}) applies.

\subsection*{${\bf n=4}$}
Choose
\begin{equation} \label{3.9}
b(t)=(1+t)^{3/5},\qquad a(t)\equiv a_0.
\end{equation}
Then
\begin{equation} \label{3.10}
R_{\bar g}(t,x)= (1+t)^{-6/5}R_h(x)\;-\;\frac{a_0^2}{4}(1+t)^{-12/5}|\Omega|_h^2(x)\;+\;\frac{6}{25}(1+t)^{-2}.
\;
\end{equation}
Using $R_h\ge 0$ and $|\Omega|_h\le \Omega_\infty$ gives
\begin{equation} \label{3.11}
R_{\bar g}(t,x)\ \ge\ \frac{6}{25}\frac{1}{(1+t)^2}\;-\;\frac{a_0^2}{4}\,\Omega_\infty^2\,\frac{1}{(1+t)^{12/5}}.
\;
\end{equation}
If $\Omega_\infty=0$ 
then $R_{\bar g}>0$ for all $t \ge 0$.
If $\Omega_\infty>0$ and
\begin{equation} \label{3.12}
\;0<a_0<\frac{2\sqrt6}{5\,\Omega_\infty}\;
\end{equation}
then $R_{\bar g}>0$ for all $t \ge 0$.

\subsection*{${\bf n=5}$}
Choose
\begin{equation} \label{3.13}
b(t)=(1+t)^{2/5},\qquad a(t)=\varepsilon(1+t)^{-2/5}.
\end{equation}
Then
\begin{equation} \label{3.14}
R_{\bar g}(t,x)= (1+t)^{-4/5}R_h(x)\;-\;\frac{\varepsilon^2}{4}(1+t)^{-12/5}|\Omega|_h^2(x)\;+\;\frac{8}{25}(1+t)^{-2}.
\;
\end{equation}
and the uniform lower bound is
\begin{equation} \label{3.15}
R_{\bar g}(t,x)\ \ge\ \frac{8}{25}\frac{1}{(1+t)^2}\;-\;\frac{\varepsilon^2}{4}\,\Omega_\infty^2\,\frac{1}{(1+t)^{12/5}}.
\;
\end{equation}
If $\Omega_\infty=0$ 
then $R_{\bar g}>0$ for all $t \ge 0$.
If $\Omega_\infty>0$ and
\begin{equation} \label{3.16}
\;0<\varepsilon<\frac{4\sqrt2}{5\,\Omega_\infty}\;
\end{equation}
then $R_{\bar g}>0$ for all $t \ge 0$.

\subsection*{${\bf n\ge 6}$}
Choose
\begin{equation} \label{3.17}
b(t)=(1+t)^{\frac{2}{n-1}},\qquad a(t)=\varepsilon(1+t)^{-\frac{n-5}{n-1}}.
\end{equation}
Then
\begin{equation} \label{3.19}
R_{\bar g}(t,x)= (1+t)^{-\frac{4}{n-1}}R_h(x)\;-\;\frac{\varepsilon^2}{4}(1+t)^{-2}|\Omega|_h^2(x)\;+\;\frac{4(n-5)}{(n-1)^2}(1+t)^{-2}.
\;
\end{equation}
Therefore the uniform lower bound is
\begin{equation} \label{3.20}
R_{\bar g}(t,x)\ \ge\ \Big(\frac{4(n-5)}{(n-1)^2}-\frac{\varepsilon^2}{4}\Omega_\infty^2\Big)\frac{1}{(1+t)^2}.
\;
\end{equation}
If $\Omega_\infty=0$ 
then $R_{\bar g}>0$ for all $t \ge 0$.
If $\Omega_\infty>0$ and
\begin{equation} \label{3.21}
\;0<\varepsilon<\frac{4\sqrt{\,n-5\,}}{(n-1)\,\Omega_\infty}\;
\end{equation}
then $R_{\bar g}>0$ for all $t \ge 0$.

\section{Discussion} \label{sect4}

With reference to Theorem \ref{thm1}, from the point of view of index theory there isn't much difference whether
the genus is one or the genus is greater than one.  In both cases $\partial \overline{M}$ is aspherical and enlargeable, and
its fundamental group satisfies the Strong Novikov Conjecture.  And in both cases there is no
index theoretic obstruction to having a Riemannian metric on
$\overline{M}$ with positive scalar curvature and positive mean curvature boundary, since such metrics exist.
Either there is no end obstruction to nonnegative scalar curvature, or index theory is not suited to find it,
or there is some new aspect to be discovered.

In general dimension, if $X$ is compact and there is a complete Riemannian metric on $[0, \infty) \times X$, with {\em uniformly} positive scalar curvature,
then $\mu$-bubbles show that there is a sequence $\{H_i\}_{i=1}^\infty$ of hypersurfaces exiting
the end so that each $H_i$ admits a Riemannian metric of positive scalar curvature, and there is a degree one
map from $H_i$ to $X$; c.f. \cite{SWWZ}. In contrast, in the proof of Theorem \ref{thm1} we obtain surfaces with a
quantitative negative lower bound on their total scalar curvatures.

If a compact manifold $X$ has a Riemannian metric of nonnegative scalar curvature then
a warped product construction shows that $[0, \infty) \times X$ has a Riemannian metric of positive scalar curvature.
However, the converse is not true, as shown by Theorem \ref{thm2} when $X$ is a three dimensional nilmanifold.
In general, the most one could expect is that $X$ has almost nonnegative scalar curvature in some sense.

To elaborate on this, a compact orientable three dimensional manifold $X$ has almost nonnegative scalar curvature, relative to the volume, if
and only if it is a graph manifold, i.e. has no hyperbolic piece in its geometric decomposition
\cite[Proposition 93.10]{Kleiner-Lott (2008)}.
As to when $[0, \infty) \times X$ has a metric of positive scalar curvature,
Theorem \ref{thm2}, or more precisely its proof, implies that this is the case when $X$ has Thurston type $\R^3$, $S^3$, $\R \times S^2$ or $\Nil$.
If $X$ has Thurston type $H^3$, $\R \times H^2$, $\widetilde{\SL(2, \R)}$ or $\Sol$ then we were unable to find such a metric.

The manifold $X$ in Theorem \ref{thm2} has almost nonnegative scalar curvature in the sense that with reference to the
metric $g$ in (\ref{3.3}),
as $a \rightarrow 0$,
i.e. as the circle fiber shrinks, the scalar curvature has a negative lower bound that approaches zero.

       \end{document}